\documentclass{article}

\usepackage[
  lang = american, 
 paper = hardcover, openaccess,
]{ems-ecm}

\numberwithin{equation}{section}

\newtheorem{thm}{Theorem}[section]

\newtheorem{defin}[thm]{Definition}

\newtheorem{conj}[thm]{Conjecture}

\theoremstyle{definition} 

\newtheorem{example}[thm]{Example}

\newtheorem{aim}[thm]{Aim}
\newtheorem{proposal}[thm]{Proposal}

\newtheorem{remark}[thm]{Remark}

\newtheorem{problem}[thm]{Problem}

\usepackage{tikz,mathrsfs}
\usetikzlibrary{arrows,decorations.pathmorphing,decorations.pathreplacing,positioning,shapes.geometric,shapes.misc,decorations.markings,decorations.fractals,calc,patterns,plotmarks}

\tikzset{
        cvertex/.style={circle,draw=black,inner sep=1pt,outer sep=3pt},
        pvertex/.style={circle,inner sep=1pt,outer sep=2pt,font=\scriptsize},
        vertex/.style={circle,fill=black,inner sep=1pt,outer sep=3pt},
        DB/.style={circle,draw=black,circle,fill=black,inner sep=0pt, minimum size=4pt},
        DBlue/.style={circle,draw=black,circle,fill=white,inner sep=0pt, minimum size=4pt},
        DW/.style={circle,draw=black,inner sep=0pt, minimum size=4pt},
        star/.style={circle,fill=yellow,inner sep=0.75pt,outer sep=0.75pt},
        tvertex/.style={inner sep=1pt,font=\scriptsize},
        gap/.style={inner sep=0.5pt,fill=white}}

\usepackage{euscript} 
\usepackage{pifont} 

\newcommand{\scrC}{\EuScript{C}}
\newcommand{\scrD}{\EuScript{D}}
\newcommand{\scrE}{\EuScript{E}}
\newcommand{\scrF}{\EuScript{F}}
\newcommand{\scrG}{\EuScript{G}}

\newcommand{\scrI}{\EuScript{I}}
\newcommand{\scrJ}{\EuScript{J}}

\newcommand{\scrL}{\EuScript{L}}

\newcommand{\scrO}{\EuScript{O}}

\newcommand{\scrR}{\EuScript{R}}

\newcommand{\scrX}{\EuScript{X}}

\newcommand\art{\mathfrak{art}}
\newcommand\cart{\mathfrak{cart}}

\newcommand\Spec{\operatorname{\mathrm{Spec}}}
\newcommand\End{\operatorname{\mathrm{End}}}
\newcommand\Hom{\operatorname{\mathrm{Hom}}}

\newcommand\Qcoh{\operatorname{\mathrm{Qcoh}}}

\renewcommand\mod{\operatorname{\mathrm{mod}}}

\newcommand\Sets{\operatorname{\mathrm{Sets}}}
\newcommand\m{\mathfrak{m}}
\newcommand\n{\mathfrak{n}}

\newcommand\cDef{\mathrm{c}\scrD\mathrm{ef}}
\newcommand\Def{\scrD\mathrm{ef}}

\newcommand{\corank}{\scrC\mathrm{rk}}

\newcommand{\dcyc}{\updelta}

\newcommand{\llangle}{\langle\kern -2.5pt\langle}
\newcommand{\rrangle}{\rangle\kern -2.5pt\rangle}
\newcommand{\llsq}{[\kern -2.5pt [}
\newcommand{\rrsq}{]\kern -2.5pt ]}
\newcommand{\lcl}{(\kern -2.5pt(}
\newcommand{\rcl}{)\kern -2.5pt)}
\newcommand{\Jac}{\scrJ\mathrm{ac}}

\newcommand{\GKdim}{\mathop{\mathrm{GKdim}}\nolimits}
\newcommand{\JRdim}{\mathop{\mathrm{Jdim}}\nolimits}
\newcommand{\frakJdim}{$\mathfrak{J}$\nobreakdash-dimension}

\numberwithin{equation}{section}

\begin{document}


\title{Noncommutative Singularity Theory}



 \emsauthor{1}{
	\givenname{Gavin}
	\surname{Brown}
	\mrid{648731} 
	\orcid{0000-0002-4087-5624}}{G.~Brown} 
\emsauthor{2}{
	\givenname{Michael}
	\surname{Wemyss}
	\mrid{893224} 
	\orcid{0000-0003-4760-7750}}{M.~Wemyss} 
 

\Emsaffil{1}{
	\department{Mathematics Institute}
	\organisation{University of Warwick}
	\rorid{01a77tt86} 
	\address{Zeeman Building}
	\zip{CV4 7AL}
	\city{Coventry}
	\country{UK}
	\affemail{G.Brown@warwick.ac.uk}}
 \Emsaffil{2}{
	\department{School of Mathematics and Statistics}
	\organisation{University of Glasgow}
	\rorid{00vtgdb53} 
	\address{University Place}
	\zip{G12 8QQ}
	\city{Glasgow}
	\country{UK}
	\affemail{michael.wemyss@glasgow.ac.uk}}

\classification[14B20, 14E30]{14B05}



\begin{abstract}
This is an expository article on the noncommutative singularity theory of power series in noncommuting variables, its motivation from deformation theory, and its applications to contractibility of curves and the classification of smooth 3-fold flops. 
\end{abstract}

\maketitle

\section{Introduction}
It is a fact of life that there are just too many functions, and to gain structural insight, some restriction is necessary. Taking that as given, what comes next is harder, since it involves making choices.  

To first approximation, making the choice of a class of functions amounts to making the choice of a research area: differential geometry is often considered to be the study of $C^\infty$ functions, algebraic geometry thought to be the exclusive realm of polynomial functions, and symplectic topology to be, well, the land where symplectomorphisms roam free.

Whilst this viewpoint provides helpful motivation, it is simply flat-out false from any reasonable technical perspective.  There are many situations where within an area `defined' by a class of functions, to prove certain theorems requires an enlargement or extension to that class of functions.  This may be frustrating. For example, if you begin with polynomials that look as nice as $u^2+v^2+x^2+y^2$, it does initially take some convincing that formal power series really are objects that humans should consider, and that analytic changes in coordinates really do make your life easier. 
If a problem can be stated entirely within one class of functions, there is understandably a reluctance to extend that class without a clear and compelling reason to do so. 

\medskip
The purpose of this article is to introduce noncommutative singularity theory as an object in its own right, and to then explain why this extension from commutative to noncommutative polynomials (and power series) is needed to solve various contraction and classification problems within algebraic geometry. 

With this in mind, and attempting to accomplish multiple goals, the structure of this article is as follows. \S\ref{subsec:NCsings} outlines what noncommutative singularity theory is and states some of its initial results. This is all factual, but it doesn't address the why: that is partially contained in  \S\ref{sec:deftheory}, which outlines the original intuition from deformation theory, but in fact is mostly explained in \S\ref{sec:applications1} and \S\ref{sec:applications2}, which covers in detail the applications to birational geometry.  Throughout \S\ref{sec:applications1} and \S\ref{sec:applications2} there is a particular focus on what the new noncommutative technology buys us, in the form of various theorems and geometric consequences. Two conjectures are briefly outlined in \S\ref{subsec:conjectures}.

\section{Noncommutative Singularity Theory}\label{subsec:NCsings}
Polynomial functions are elements in the polynomial ring $\mathbb{C}[ x_1,\hdots,x_d]$ for some $d\geq 1$. To set language, noncommutative (polynomial) functions will simply be elements of the free algebra $\mathbb{C}\langle x_1,\hdots,x_d\rangle$.

\subsection{Classical singularity theory}\label{subsec:classicalsing}
To classify even smooth points in algebraic geometry, one must work suitably locally.  
For this article this means working with elements of the commutative power series ring $\mathbb{C}\llsq x_1,\hdots,x_d\rrsq$, although alternatives are available.  The key point is that elements in the commutative power series ring are allowed to consist of infinite sums, and as such they can be represented by the form
\[
f=\uplambda_1 + \uplambda_2x+ \uplambda_3 y + \uplambda_4 xy +\hdots
\]
The ring $\mathbb{C}\llsq x_1,\hdots,x_d\rrsq$ is sometimes alternatively referred to as the formal power series ring, with the word `formal' acting as a reminder that its elements do not need to satisfy any radius of convergence conditions. 

At first glance, it is not so clear that moving from polynomials (finite sums) to power series (infinite sums) makes life easier, given that there seem to be many more elements.  But it does, since there are now \emph{many} more symmetries, in the form of analytic changes in coordinates.  Given this huge symmetry group, it is then possible to ask about classifying functions.  This is the domain of singularity theory.

\begin{aim}[Singularity Theory \`a la Arnold \cite{Arnold}]\label{aim:Arnold}
Up to specified isomorphism, classify all $f\in\mathbb{C}\llsq x_1,\hdots,x_d\rrsq$ that satisfy some fixed numerical criteria, and produce theory for when classification is not possible.
\end{aim}

The following is perhaps the most fundamental class of examples.   Although the condition below is a priori unmotivated, it is numerical, being the count of ideals satisfying a certain property. 
\begin{defin}\label{def:simple}
An element $f\in\mathbb{C}\llsq x_1,\hdots,x_d\rrsq=\mathbb{C}\llsq \mathsf{x}\rrsq$ is said to be a simple singularity if
\[
\# \{ I \mid I \textnormal{ proper ideal of } \mathbb{C}\llsq \mathsf{x}\rrsq \textnormal{ with }f \in I^{2} \}<\infty.
\]
\end{defin}
This is a non-standard definition of a simple singularity.  Like many numerical definitions within singularity theory, it is only really justified \emph{after} classification, when it becomes more clear that it is equivalent to every other characterisation of simple singularities in the literature.
\begin{thm}[{e.g.\ \cite[\S8]{Yoshino}}]\label{thm:ADE}
If $f$ is a simple singularity in the sense of \ref{def:simple}, then up to relabelling variables $z_1,\hdots z_{d-2},x,y$, and up to isomorphism (analytic or formal changes of coordinates in the domain), $f$ is one of
\[
\begin{array}{lll}
 A_n& \quad z^{2} + x^{2} + y^{n+1}& n\geq 1 \\[1mm]
 D_n&\quad z^{2} + x^{2}y + y^{n-1}& n\geq 4\\[1mm]
 E_6&\quad z^{2} + x^{3} + y^{4}& \\[1mm]
 E_7&\quad z^{2} + x^{3} + xy^{3}&  \\[1mm]
 E_8&\quad z^{2} + x^{3} + y^{5}&  
\end{array}
\]
where $z^{2}=z_1^{2}+\hdots+ z_{d-2}^{2}$.
\end{thm}
From the viewpoint of later sections, several remarks are in order. First, in \ref{def:simple} there is  no mention of ADE, yet this still appears in the classification \ref{thm:ADE}.  Second, \ref{thm:ADE} illustrates the very small number of simple singularities, namely just two (countable) infinite families and only three other cases. And third, with \ref{thm:ADE} in hand, working with simple singularities becomes much easier, since various problems may be reduced to checking statements against a short list of elegant forms.

\subsection{Noncommutative singularity theory} 
To first approximation, the overall idea of noncommutative singularity theory is to replace the commutative object $\mathbb{C}\llsq x_1,\hdots,x_d\rrsq$ with the formal \emph{noncommutative} power series ring $\mathbb{C}\llangle x_1,\hdots,x_d\rrangle$, then to play the same game as in \S\ref{subsec:classicalsing}.  This is not quite true, since in the noncommutative setting there are many more `free objects' given by quivers, and as in \ref{rem:quivers} below this makes the theory richer, and widens applications.

 In the first case, the ring $\mathbb{C}\llangle x_1,\hdots,x_d\rrangle$ is similar to the free algebra in $d$ variables, except now elements are possibly infinite sums
\[
f=\uplambda_1 + \uplambda_2x+ \uplambda_3 y + \uplambda_4 xy +\uplambda_5 yx + \hdots
\]
Whilst for $d\geq 2$ the ring $\mathbb{C}\llangle x_1,\hdots,x_d\rrangle$ is not noetherian, has exponential growth, and is pathologically horrible from many perspectives, as in the commutative case its automorphim group is huge, and so classifying elements satisfying certain numerical criteria still turns out to be possible. But, as in the classical case, to do this still requires choices.

\begin{problem}
When should $f, g\in\mathbb{C}\llangle x_1,\hdots,x_d\rrangle$ be viewed as being the same? And what would be natural numerical criteria?
\end{problem}

To explain first the notion of equivalence, recall that every $f\in\mathbb{C}\llangle x_1,\hdots,x_d\rrangle$ can be differentiated with respect to any variable.  The rule is simple: cycle then score off the leftmost variable.  For example, in humane coordinates $x$ and $y$, to differentiate $x^3y$ with respect to $x$,
\[
x^{3}y=xxxy\xrightarrow{\scriptstyle \textnormal{cycle}}
\begin{array}{c}xxxy\\xxyx \\xyxx\\yxxx
\end{array}
\xrightarrow{\scriptstyle \textnormal{score off}}
\begin{array}{c}xxy\\xyx \\yxx\\0
\end{array}
\]
where `score off' either scores off the leftmost $x$, or returns zero if the monomial does not begin with $x$.  By definition, it follows that  $\updelta_x(x^{3}y)=x^{2}y+xyx+yx^{2}$.  Reassuringly, $\updelta_y(y^{3})=3y^{2}$.

\begin{defin}
Given any $f$, the Jacobi algebra is
\[
\Jac f:= \frac{\mathbb{C}\llangle x_1,\hdots,x_d\rrangle}{\lcl \updelta_zf \mid z \mbox{ is a variable}\rcl},
\]
where the double brackets $\lcl\,\rcl$ denote the closure of the two-sided ideal.
\end{defin}

\begin{example}\label{example9dim}
The element $f=x^4 + xy^2\in \mathbb{C}\llangle x,y\rrangle$ has
\[
\Jac f=\frac{\mathbb{C}\llangle x,y\rrangle}{\lcl 4x^3 + y^2, xy+yx\rcl},
\]
which as explained in \ref{ex:9dB} below is a nine-dimensional algebra.
\end{example}

Write $f\cong g$ provided that $\Jac f\cong\Jac\, g$.  

\begin{aim}
With $\mathbb{C}\llangle x_1,\hdots,x_d\rrangle$ fixed, and the notion of equivalence $\cong$ now defined, the aim of singularity theory remains: classify all elements satisfying some numerical (or structural) criteria, and produce theory for when classification is not possible.
\end{aim}

\begin{remark}\label{rem:quivers}
An alternative description of the noncommutative power series ring $\mathbb{C}\llangle x_1,\hdots,x_d\rrangle$ is the complete path algebra of the $d$-loop quiver.  Other quivers are possible, and all of the constructions in this section and below naturally and easily generalise to that setting \cite{DWZ}.  These other quivers are needed for multi-curve geometric applications; see \S\ref{NC multiple sction}, \S\ref{deform curves} and \ref{rem:Hao} later.
\end{remark}

\subsection{The singularity theory toolkit}
This subsection briefly outlines some of the singularity theory toolkit, adapted into the noncommutative setting.

The first most fundamental problem is to establish the correct notion of dimension and the related numerical invariants.  Using the Gelfand--Kirillov (GK) dimension, which is standard for noncommutative rings, is problematic since $\GKdim\mathbb{C}\llsq x\rrsq=\infty$.  Fundamentally, GK dimension is not well adapted to complete rings. To overcome this, we return to the basic fact that $\Jac f$ is local, with Jacobson radical $\mathfrak{J}$ say.  This gives an intrinsic filtration
\[
\Jac f\supseteq\mathfrak{J}\supseteq\mathfrak{J}^2\supseteq\hdots.
\] 
\begin{defin}
The (Jacobson-radical) \frakJdim, $\JRdim\Jac f$, is the growth rate of the steps in the chain above, namely
\begin{equation}
\inf \left\{ r\in\mathbb{R} \mid 
\textnormal{for some $c\in\mathbb{R}$, }\dim\Jac f / \mathfrak{J}^n \le cn^r \text{ for every $n\in\mathbb{N}$} \right\}.\label{eqn:Jacobchain}
\end{equation}
\end{defin}
To calibrate, note that $\JRdim\mathbb{C}\llsq x\rrsq=1$, and that $\JRdim\Jac f=0$ is equivalent to $\Jac f$ being a finite dimensional algebra.

\begin{remark}
There are several well-known results about GK dimension: for example, that it can take any real value $\ge2$ (or $0$, $1$ or $\infty$) but no values in the `Bergman gap' between $1$ and~$2$. We do not know analogous results for \frakJdim, even for the restricted class of Jacobi algebras. In the context of this survey, however, it is clear that if the \frakJdim\, is $\le1$ then it is either~0 or~1, just as for GK dimension.
\end{remark}

Whilst the growth rate of the whole chain \eqref{eqn:Jacobchain} gives the dimension, the numerical properties of its individual pieces also holds valuable information.  In what follows, we will always assume that $f$ contains no constant or linear terms, which we will abbreviate as $f\in\mathbb{C}\llangle x_1,\hdots,x_d\rrangle_{\geq 2}$.

\begin{defin}\label{defin corank}
For $f\in\mathbb{C}\llangle x_1,\hdots,x_d\rrangle_{\geq 2}$, and $m\ge1$, the {\em $m$th corank of $f$} is defined to be
\[
\corank_m(f) =\dim_{\mathbb{C}}\left(\frac{\mathfrak{J}^m}{\mathfrak{J}^{m+1}}\right),
\]
where $\mathfrak{J}$ is the Jacobson radical of $\Jac f$. 
\end{defin}
The first of these invariants $\corank(f):=\corank_1(f)$ is often referred to as simply the corank.  The corank is determined by the linear conditions imposed by derivatives, and is therefore uniquely determined by the degree two part of $f$.

\begin{thm}[{Splitting Lemma, \cite[4.5]{DWZ}}]
Let $f\in\mathbb{C}\llangle x_1,\hdots,x_d\rrangle_{\geq 2}$. Then $f\cong x_1^2+\cdots+x_r^2+g$ for some $g\in\mathbb{C}\llangle x_{r+1},\hdots, x_d\rrangle_{\ge3}$, where $d-r=\corank(f)$. In particular,
\[
\Jac f \cong \frac{\mathbb{C}\llangle x_{r+1},\hdots, x_d\rrangle}{\lcl \dcyc_{x_{r+1}}g,\hdots,\dcyc_{x_d}g \rcl}.
\]
\end{thm}

The coranks are properties of the leading terms of elements of the closed Jacobian ideal of derivatives, and can be calculated using a standard basis. In the commutative context, the standard basis (also known as a local Gr\"obner basis) is well known. Loosely speaking, one follows the usual Buchberger approach of cancelling leading terms, but regarding the leading terms as those of low degree. The theory and key results all hold in the setting of noncommutative power series \cite{GH98}. From a user persepective, there are two points to note: noncommutativity means leading terms may have nontrivial cancellation with themselves, and, being non-noetherian, the resulting basis may be infinite (although not of course in nice settings such as $\JRdim \Jac f=0$).

\begin{example}\label{ex:9dB}
    Continuing \ref{example9dim}, consider $f=x^4+xy^2$. The leading term of a nonzero power series is defined to be the term of least degree, with lexicographic order used to break ties. So, for example, the leading term of $xy^2+yxy+y^3$ is $xy^2$.

    Setting $A=y^2+4x^3$ and $B=xy+yx$ for the two derivatives, the leading term $y^2$ of $A$ cancels itself in the expression $yA-Ay$ to give
    \[
    yA-Ay = 4yx^3 - 4x^3y \buildrel{B}\over{\equiv} 4yx^3 + 4x^2yx \buildrel{B}\over{\equiv} 4yx^3 - 4xyx^2 \buildrel{B}\over{\equiv} 4yx^3 + 4yx^3 = 8yx^3
    \]
    where $\buildrel{B}\over{\equiv}$ denotes congruence modulo~$B$. Thus $C=yx^3\in (A,B)$. Similarly $Ax^3-yC=4x^6\in(A,B)$. All other attempts to cancel leading terms reduce to zero, and so the (normalised) standard basis is
    \[
    y^2 + 4x^3,\ xy+yx,\ yx^3,\ x^6.
    \]
    The leading terms of these are $y^2,xy,yx^3,x^6$, and the set of monomials not divisible by any of these provides a basis of $\Jac f$. In this case we see at once therefore that
    \[
    1,\ x,\ y,\ x^2,\ yx,x^3,\ yx^2,\ x^4,\ x^5
    \]
    is a basis of $\Jac f $ and so it has dimension~9, as claimed.
    Morover, the degree~$m$ monomials in this monomial basis descend to a basis of $\mathfrak{J}^m/\mathfrak{J}^{m+1}$, which implies that
    \[
    \corank_1(f) = \corank_2(f)=\corank_3(f)=2 \textrm{\ \ and\ \ } 
    \corank_4(f)=\corank_5(f)=1.
    \]
    Together with $\dim_{\mathbb{C}} \mathbb{C}\llangle x, y\rrangle/\mathfrak{J}=1$, these numbers add up to~$9=\dim_{\mathbb{C}}\Jac f$.
\end{example}

Another key property of Jacobi algebras is that they are invariant under automorphisms, and so in particular under arbitrary analytic changes in coordinates.
\begin{thm}[{\cite[3.7]{DWZ}}]\label{thm:autos}
If $F\in\mathrm{Aut}\,\mathbb{C}\llangle x_1,\hdots,x_d\rrangle $ then $\Jac f\cong \Jac\, F(f)$.
\end{thm}

Other natural parts of the commutative theory, such as Saito's famous result \cite{S} on quasihomogeneous power series, also work in the noncommutative context, although not quite yet at the elementary entry level of this discussion; see \cite{HZ}. 

\begin{remark}
Whilst any type of singularity theory is really `just' the classification of Jacobi (or commutatively, Tjurina) algebras up to isomorphism, thinking in terms of functions allows us to access the toolkit of singularity theory.  Even when the Jacobi algebras are finite dimensional, perhaps the main advantage of thinking in terms of functions is \ref{thm:autos}, as this enables the construction of many isomorphisms which seem completely out of reach of more classical approaches to the classification of noncommutative finite dimensional  algebras.
\end{remark}

\subsection{Proposal}
The following proposal \ref{prop:NCsing}, which we view as the noncommutative version of simple singularities, is our main reason for developing noncommutative singularity theory.

\begin{proposal}\label{prop:NCsing}
Under the numerical criterion $\JRdim\Jac f\leq 1$, we propose that:
\begin{enumerate}
\item a classification is in fact possible.
\item the classification is ADE.
\item this algebraic classification \emph{is} (and implies) the classification of flops and of crepant divisorial contractions to curves.
\end{enumerate}
\end{proposal}

The proposal is still conjectural in its form stated above, as this is strongly related to Conjecture~\ref{conj:georeal} below. The current state of the art, partially explained in \S\ref{sec:applications2}, in part replaces the numerical condition $\JRdim\Jac f\leq 1$ with a structural one, and solves that: for all geometric applications including the classification of $3$-fold flops, this is enough. However, from the viewpoint of future ease of use (numerical properties are easier to verify than structural ones), and also to make the analogy with simple singularities in \ref{thm:ADE} more compelling, verifying \ref{prop:NCsing} in the generality stated remains the founding question of noncommutative singularity theory.

\section{Deformation Theory}\label{sec:deftheory}

This section describes the links, via deformation theory, between noncommutative singularity theory on the one hand and both classification and contractibility problems in birational geometry on the other.  Applications are given later, in \S\ref{sec:applications1} and \S\ref{sec:applications2}.

\subsection{Deformation theory framework}
One reason for deforming an object, say $M$, is the hope that doing so will attach an intrinsic invariant to $M$, and furthermore that this invariant will be useful.  In good cases it is possible to extract simple data from the invariant, for example a number such as its dimension (if that even makes sense). Intrinsically attaching numbers to objects is hard, and the framework of deformation theory is generally a good way of doing this.

Alas, deformation theory involves choices.  To ensure that the theory is intrinsic, it is natural to insist on the functorial approach.  That is, given some object $M$, the first aim is to construct a functor
\[
\scrD\mathrm{ef}_M\colon \mathrm{TestObjects}\to\mathrm{Sets}
\]
which sends a test object $R$ to a set $\{\mbox{deformations of }M\mbox{ over }R\}/\sim$. To make this rigorous, the following need to be defined, or rather declared:
\begin{enumerate}
\item The category $\mathrm{TestObjects}$.
\item The definition of `deformations of $M$ over $R$'.
\item The equivalence relation $\sim$
\end{enumerate}
There is a fourth choice, suppressed throughout this article, since we will always assume that the target category is the category of sets.

\subsection{Commutative deformation theory}\label{subsec:commdef}
Let $X$ be a variety, and $\scrF\in\Qcoh X$ be a quasi-coherent sheaf.  To deform $\scrF$ via the classical Grothendieck formulation, consider the functor 
\[
\cDef\colon \cart\to\Sets
\]
where $\cart$ is the category of commutative local artinian $\mathbb{C}$-algebras, sending
\[
(R,\m)\mapsto \left. \left \{ 
(\scrG,\upeta)
\left|\begin{array}{l}
\scrG\in\Qcoh (X\times \Spec R),\\ 
\scrG\mbox{ is flat over }\Spec R,\\ 
\upeta\colon \scrG|_{X\times 0}\xrightarrow{\sim}\scrF.
\end{array}\right. \right\} \middle/ \sim \right.
\]
where $\sim$ is the natural equivalence relation outlined in a greater generality below.  The point is that, roughly speaking, a deformation of $\scrF$ is a sheaf $\scrG$ which is infinitesimally thicker than $\scrF$ (namely, a sheaf on $X\times \Spec R$), roughly of the same size (flatness), which recovers the original $\scrF$ (restricted to the zero fibre).  There are many choices of such isomorphisms on the zero fibre, so we carry those choices as part of the data.

\subsection{Noncommutative deformation theory}\label{subsec:NCdef}
It is very intuitive that more test objects will lead to more information.  This subsection enlarges the category of test objects $\cart$ from \S\ref{subsec:commdef} by adding in (some) noncommutative artinian rings.  The resulting deformation functor will be similar, but requires some reinterpretation.

For simplicity, consider the same setup as in \S\ref{subsec:commdef}, namely $X$ is a variety, and $\scrF\in\Qcoh X$.  The noncommutative deformation functor is defined to be
\[
\Def\colon \art_1\to\Sets
\]
from finite dimensional augmented $\mathbb{C}$-algebras to sets, sending
\[
(\Gamma,\n)\mapsto \left. \left \{ 
(\scrG,\upvartheta,\upeta)
\left|\begin{array}{l}
\scrG\in\Qcoh X,\\
\upvartheta\colon \Gamma\to\End_{X}(\scrG)\mbox{ homomorphism},\\ 
-\otimes_\Gamma \scrG\colon\mod\Gamma\to \Qcoh X \mbox{ is exact},\\ 
\upeta\colon (\Gamma/\n)\otimes_\Gamma\scrG\xrightarrow{\sim}\scrF.
\end{array}\right. \right\} \middle/ \sim \right.
\]
where, for any given $\Gamma$, $(\scrG,\upvartheta,\upeta)\sim (\scrG',\upvartheta',\upeta')$ if and only if there exists $\uptau \colon \scrG\xrightarrow{\sim}\scrG'$ such that for all $r\in\Gamma$
\[
\begin{array}{c}
\begin{tikzpicture}
\node (a1) at (0,0) {$\scrG$};
\node (a2) at (1.5,0) {$\scrG'$};
\node (b1) at (0,-1.5) {$\scrG$};
\node (b2) at (1.5,-1.5) {$\scrG'$};
\draw[->] (a1) -- node[above] {$\scriptstyle \uptau$} (a2);
\draw[->] (a1) -- node[left] {$\scriptstyle \upvartheta(r)$} (b1);
\draw[->] (a2) -- node[right] {$\scriptstyle \upvartheta'(r)$} (b2);
\draw[->] (b1) -- node[above] {$\scriptstyle \uptau$} (b2);
\end{tikzpicture}
\end{array}
\quad
\mbox{and}
\quad
\begin{array}{c}
\begin{tikzpicture}
\node (a1) at (0,0) {$(\Gamma/\n)\otimes_\Gamma \scrG$};
\node (a2) at (3,0) {$(\Gamma/\n)\otimes_\Gamma \scrG'$};
\node (b) at (1.5,-1.5) {$\scrF$};
\draw[->] (a1) -- node[above] {$\scriptstyle 1\otimes \kern 1pt \uptau$} (a2);
\draw[->] (a1) -- node[pos=0.6,left,inner sep=8,anchor=east] {$\scriptstyle \upeta$} (b);
\draw[->] (a2) -- node[pos=0.6,right,inner sep=8,anchor=west] {$\scriptstyle \upeta'$} (b);
\end{tikzpicture}
\end{array}
\]
commute.  There are many equivalent characterisations of the category $\art_1$, where we refer the reader to \cite{Eriksen2} for full details.

The commutative deformation functor is recovered via 
\[
\cDef=\Def|_{\cart}.
\]

\subsection{Deforming multiple objects}\label{NC multiple sction}
Perhaps the clearest way of observing the benefit of noncommutative deformation theory is when we try to deform more than one object, say $\scrF_1$ and $\scrF_2$.  Classically, the commutative solution is to take the direct sum $\scrF_1\oplus \scrF_2$, and apply the technology of \S\ref{subsec:commdef}. The commutative hull is then constructed as a factor of $\mathbb{C}\llsq x_1,\hdots,x_t\rrsq$, where $t=\dim_{\mathbb{C}}\mathrm{Ext}^1_X(\scrF_1\oplus \scrF_2,\scrF_1\oplus \scrF_2)$.

But this clearly loses information, since Ext groups are directed: indeed, in general
\[
\mathrm{Ext}^1_X(\scrF_1,\scrF_2)\ncong \mathrm{Ext}^1_X(\scrF_2,\scrF_1).
\]
Since the order matters in the Ext groups, forcing the commutativity in $\mathbb{C}\llsq x_1,\hdots,x_t\rrsq$ means that often almost all information is lost.

The solution to deforming $n$ objects is to realise that, since Exts are directed, the prorepresenting hull should naturally be constructed as a factor of a quiver with $n$ vertices, where the number of arrows from vertex $i$ to vertex $j$ should be the integer $\dim_{\mathbb{C}}\mathrm{Ext}^1_X(\scrF_i,\scrF_j)$. This naturally records the asymmetry in the Ext groups, and also gives a more natural and compelling description of what the tangent space should be.

To do this on a technical level requires the introduction of the category of test objects $\art_n$, which are roughly speaking finite dimensional $\mathbb{C}$-algebras where all simples are one-dimensional; this last condition rules out matrix algebras appearing.  For many equivalent characterisations of the category $\art_n$, we refer the reader to \cite{Eriksen2}.

\medskip
To deform the collection $\scrF_1,\hdots,\scrF_n$, each $\scrF_i$ is replaced by an injective resolution $\scrF_i\hookrightarrow\scrI_i$, and from that the endomorphism DG-algebra $\scrE\mathrm{nd}_X(\oplus\scrI_i)$ is constructed.  The Maurer--Cartan formulation then gives a functor
\begin{equation}
\Def_{\scrF_1,\hdots,\scrF_n}\colon\art_n\to\Sets\label{eqn:NCdefkey}
\end{equation}
called the multi-pointed noncommutative deformation functor.  Full details can be found in e.g.\ \cite{Segal, DW2}.  As always there are choices; it is also possible to consider 
\[
\Def_{\scrF_1\oplus\hdots\oplus\scrF_n}\colon\art_1\to\Sets
\]
from \S\ref{subsec:NCdef}. We describe in \S\ref{sec:applications1} below the behaviour of these functors in the situation of curve contractions, and outline the reasons why from almost all geometric viewpoints \eqref{eqn:NCdefkey} is the key.

\section{The setting of curves}\label{sec:applications1}

Our main interest is in formal or analytic neighbourhoods of rational curves on $3$-folds.

\subsection{Contractions}\label{deform curves}
Pictorially, we will consider contractions of $3$-folds, of two types:
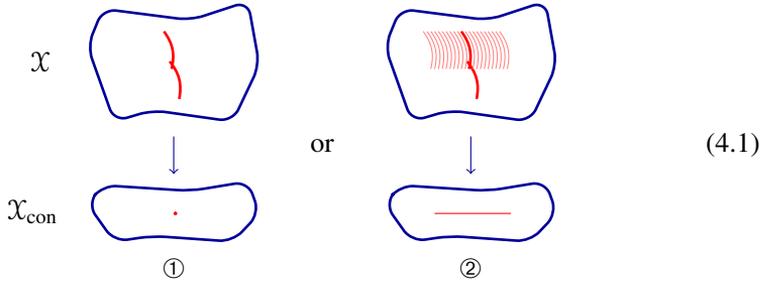
\begin{equation}
\begin{array}{ccc}
\begin{array}{c}
\begin{tikzpicture}
\node at (-1.75,0) {$\scrX$};
\node at (-1.85,-2) {$\scrX_{\mathrm{con}}$};
\node at (0,0) {\begin{tikzpicture}[scale=0.5]
\coordinate (T) at (1.8,2.4);
\coordinate (B) at (2,1.4);
\coordinate (Bb) at (1.95,1.6);
\coordinate (A) at (2.2,0.6);
\draw[red,line width=1pt] (T) to [bend left=25] (B);
\draw[red,line width=1pt] (Bb) to [bend left=25] (A);
\draw[color=blue!60!black,rounded corners=5pt,line width=1pt] (0.5,0) -- (1.5,0.3)-- (3.6,0) -- (4.3,1.5)-- (4,3.2) -- (2.5,2.7) -- (0.2,3) -- (-0.2,2)-- cycle;
\end{tikzpicture}};
\node at (0,-2) {\begin{tikzpicture}[scale=0.5]
\filldraw [red] (2.1,0.75) circle (1pt);
\draw[color=blue!60!black,rounded corners=5pt,line width=1pt] (0.5,0) -- (1.5,0.15)-- (3.6,0) -- (4.3,0.75)-- (4,1.6) -- (2.5,1.35) -- (0.2,1.5) -- (-0.2,1)-- cycle;
\end{tikzpicture}};
\draw[->, color=blue!60!black] (0,-1) -- (0,-1.5);
\node at (0,-2.75) {\ding{192}};
\end{tikzpicture}
\end{array}
& \mbox{or} &
\begin{array}{c}
\begin{tikzpicture}
\node at (0,0) {\begin{tikzpicture}[scale=0.5]
\coordinate (T) at (1.8,2.4);
\coordinate (B) at (2,1.4);
\coordinate (Bb) at (1.95,1.6);
\coordinate (A) at (2.2,0.6);
\draw[red,line width=1pt] (Bb) to [bend left=25] (A);

\draw[red,line width=1pt] (T) to [bend left=25] (B);
\foreach \y in {0.1,0.2,...,1}{ 
\draw[very thin,red!50] ($(T)+(\y,0)$) to [bend left=25] ($(B)+(\y,0)$);
\draw[very thin,red!50] ($(T)+(-\y,0)$) to [bend left=25] ($(B)+(-\y,0)$);}
\draw[color=blue!60!black,rounded corners=5pt,line width=1pt] (0.5,0) -- (1.5,0.3)-- (3.6,0) -- (4.3,1.5)-- (4,3.2) -- (2.5,2.7) -- (0.2,3) -- (-0.2,2)-- cycle;
\end{tikzpicture}};
\node at (0,-2) {\begin{tikzpicture}[scale=0.5]
\draw [red] (1.1,0.75) -- (3.1,0.75);
\draw[color=blue!60!black,rounded corners=5pt,line width=1pt] (0.5,0) -- (1.5,0.15)-- (3.6,0) -- (4.3,0.75)-- (4,1.6) -- (2.5,1.35) -- (0.2,1.5) -- (-0.2,1)-- cycle;
\end{tikzpicture}};
\draw[->, color=blue!60!black] (0,-1) -- (0,-1.5);
\node at (0,-2.75) {\ding{193}};
\end{tikzpicture}
\end{array}
\end{array}\label{eqn:twotypes}
\end{equation}
where the contraction map $f$ is projective birational, satisfying $\mathrm{R}f_*\scrO_{\scrX}=\scrO_{\scrX_{\mathrm{con}}}$.  In the left hand case \ding{192}, the morphism is an isomorphic away from a point, whereas in the right hand case \ding{193} the morphism contracts a divisor to a curve, and is an isomorphism away from the one-dimensional locus on the base.  Although the above pictures are drawn in the algebraic category, we will often implicitly pick a closed point on the base within the non-isomorphism locus, complete at that point, and thus assume that the base $\scrX_{\mathrm{con}}=\mathrm{Spec}\,\scrR$ is complete local.

\medskip
Above the chosen point $p$, the fibre (given reduced scheme structure) is a collection of curves $C_1,\hdots,C_n$, all isomorphic to $\mathbb{P}^1$.  The number of curves depends on which closed point on the base is chosen.

Regardless, applying \S\ref{NC multiple sction} to the collection of structure sheaves $\scrO_{C_1},\hdots,\scrO_{C_n}\in\Qcoh\scrX$ gives rise to the noncommutative deformation functor 
\[
\Def_{C_1,\hdots,C_n}\colon\art_n\to\Sets.
\]

\begin{thm}[\cite{DW1,DW2}]\label{Acondefthm}
There is an algebra $\mathrm{A}_{\mathrm{con}}$ and a functorial isomorphism
\[
\Def_{C_1,\hdots,C_n}(\Gamma)\cong \Hom(\mathrm{A}_{\mathrm{con}},\Gamma)/\sim
\]
where $\sim$ is the action of inner automorphisms in $\Gamma$.
\end{thm}
The algebra $\mathrm{A}_{\mathrm{con}}$ is called the contraction algebra. In fact, there is a second algebra $\mathrm{A}_{\mathrm{def}}$ called the noncommutative deformation algebra, but it is shown in \cite{DW1} that the two are isomorphic; from the viewpoint of \S\ref{subsec:conjectures} below this latter algebra $\mathrm{A}_{\mathrm{def}}$ is important, since it can be defined \emph{without} needing to assume the existence of a contraction.  Regardless, in general $\mathrm{A}_{\mathrm{con}}$ is very challenging to compute.  At least with some extra smoothness assumptions, \S\ref{subsec:linktoNC} explains the link to noncommutative singularity theory.

We briefly note that there is a small technical issue in \cite{DW1,DW3} since unlike \cite[2.11]{Segal} those papers forgot to account for the action of inner automorphisms in the statement of \ref{Acondefthm}. In practice this does not effect anything; as explained in \cite[2.14]{Segal} even after accounting for inner automorphisms the prorepresenting object is still unique up to non-unique isomorphism.  An alternative approach to this is to redefine the morphisms in the category $\art_n$ to be conjugacy classes, as in \cite{BB}.

\subsection{Which deformation theory?}
Given the abundance of deformation functors in \S\ref{sec:deftheory}, it seems natural to wonder which solves the geometric problems of interest.

Around the chosen point $p$ on the base, above which there are curves $C_1,\hdots,C_n$, we first ask which deformation functor detects the difference between sitation \ding{192} and \ding{193}. In all cases, without assumptions on the singularities of $\scrX$, and without crepant assumptions, this is answered by the noncommutative deformation functor.

\begin{thm}[\cite{DW2}] Given a contraction $\scrX\to\scrX_{\mathrm{con}}$ as in \eqref{eqn:twotypes}, pick a point $p$ on the base and consider the corresponding contraction algebra $\mathrm{A}_{\mathrm{con}}$ above that point.  Then, locally around $p$, 
\[
\begin{array}{rcl}
\JRdim\mathrm{A}_{\mathrm{con}}=0 &\iff& \mbox{Situation \ding{192} (flopping)}.\\[2mm]
\JRdim\mathrm{A}_{\mathrm{con}}=1 &\iff& \mbox{Situation \ding{193} (div$\to$curve)}.\\
\end{array}
\]
\end{thm}
With extra smoothness assumptions, $\mathrm{A}_{\mathrm{con}}$ is a Jaocbi algebra (see \S\ref{subsec:linktoNC}), and so in addition to detecting the nature of the geometry locally around $p$, the above theorem strongly motivates the study of potentials $f$ satisfying the numerical condition $\JRdim\Jac f\leq 1$.

The larger question is which deformation functor, if any, is sufficient to \emph{classify} the curve neighbourhood, since this would be the strongest property of all. The following has become known as the Donovan--Wemyss conjecture, where for ease of exposition here we restrict to the case of a single curve.  Classification requires some form of smoothness, and the rough slogan is that contraction algebras are always the classifying objects for smooth crepant contractible neighbourhoods.

\begin{conj}{(Donovan--Wemyss)}\label{thm:DWconj}
Let $\scrX_1\to\Spec \scrR_1$  and $\scrX_2\to\Spec \scrR_2$ be as in \eqref{eqn:twotypes}, where both are crepant, both $\scrX_i$ are smooth, and both $\scrR_i$ are complete local cDV.  Write $\mathrm{A}_{\mathrm{con}}$ and $\mathrm{B}_{\mathrm{con}}$ for their corresponding contraction algebras, and in each case assume that there is only one curve.  Then
\[
\scrR_1\cong \scrR_2\iff \mathrm{A}_{\mathrm{con}}\cong\mathrm{B}_{\mathrm{con}}.
\]
\end{conj}

In the setting of flops \ding{192}, Conjecture~\ref{thm:DWconj} is now a theorem due to \cite{DW1, August, HuaKeller, JM}, where after a series of reduction steps \cite{JM} apply their Derived Auslander--Iyama correspondence.  The divisor-to-curve case of \ref{thm:DWconj} remains open.  

Regardless, it is useful to regard the word conjecture as a challenge to extend the scope of \ref{thm:DWconj} to the maximum extent possible; it may even be true that contractibility is not required, and that the noncommutative deformation algebra classifies smooth crepant curve neighbourhoods in full generality.

\subsection{Link to noncommutative singularity theory}\label{subsec:linktoNC}
As already alluded to, when $\scrX$ is smooth and the morphism is crepant, it is automatic that $\mathrm{A}_{\mathrm{con}}$ is the Jacobi algebra of a quiver with potential.  The quiver can be described easily from the data of the configuration of the curves contracted and their normal bundles \cite[2.15]{HomMMP}. Describing the potential is significantly more challenging.  

This is precisely the domain of noncommutative singularity theory: with the combinatorial data of the quiver fixed, it asks what functions (=potentials) can arise.  To calibrate, in the case of a single crepant curve, which is often the hardest case, $\mathrm{A}_{\mathrm{con}}$ is a factor of $\mathbb{C}\llangle x,y\rrangle$.  The fact we should expect only two variables to show up should be viewed in analogy with simple singularities in \ref{thm:ADE}; for a full explanation, see \cite[A.13]{BW2}.

\section{New Applications and Results}\label{sec:applications2}

This section briefly summarises some of the upcoming work \cite{BW3} which applies noncommutative singularity theory to classify normal forms into noncommutative analogues of the ADE families. We focus here on the properties that we currently know to be numerical, leaving the technical aspects to \cite{BW3}.  In \S\ref{subsec:geocor} we then briefly outline the main application to the classification of smooth $3$-fold flops, and in \S\ref{subsec:curveinv} outline some of the more surprising consequences to curve counting invariants.

To simplify the discussion, henceforth we consider crepant contractions of $3$-folds of only the following two types:
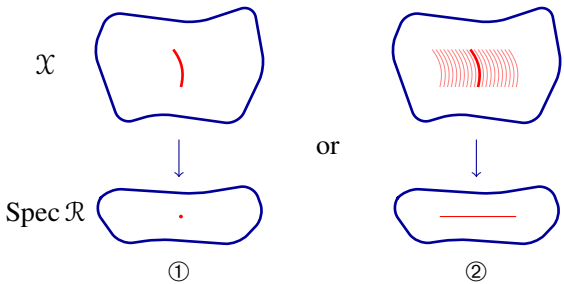
\[
\begin{array}{ccc}
\begin{array}{c}
\begin{tikzpicture}
\node at (-1.75,0) {$\scrX$};
\node at (-1.75,-2) {$\mathrm{Spec}\,\scrR$};
\node at (0,0) {\begin{tikzpicture}[scale=0.5]
\coordinate (T) at (1.9,2);
\coordinate (B) at (2.1,1);
\draw[red,line width=1pt] (T) to [bend left=25] (B);
\draw[color=blue!60!black,rounded corners=5pt,line width=1pt] (0.5,0) -- (1.5,0.3)-- (3.6,0) -- (4.3,1.5)-- (4,3.2) -- (2.5,2.7) -- (0.2,3) -- (-0.2,2)-- cycle;
\end{tikzpicture}};
\node at (0,-2) {\begin{tikzpicture}[scale=0.5]
\filldraw [red] (2.1,0.75) circle (1pt);
\draw[color=blue!60!black,rounded corners=5pt,line width=1pt] (0.5,0) -- (1.5,0.15)-- (3.6,0) -- (4.3,0.75)-- (4,1.6) -- (2.5,1.35) -- (0.2,1.5) -- (-0.2,1)-- cycle;
\end{tikzpicture}};
\draw[->, color=blue!60!black] (0,-1) -- (0,-1.5);
\node at (0,-2.75) {\ding{192}};
\end{tikzpicture}
\end{array}
& \mbox{or} &
\begin{array}{c}
\begin{tikzpicture}
\node at (0,0) {\begin{tikzpicture}[scale=0.5]
\coordinate (T) at (1.9,2);
\coordinate (B) at (2.1,1);
\draw[red,line width=1pt] (T) to [bend left=25] (B);
\foreach \y in {0.1,0.2,...,1}{ 
\draw[very thin,red!50] ($(T)+(\y,0)$) to [bend left=25] ($(B)+(\y,0)$);
\draw[very thin,red!50] ($(T)+(-\y,0)$) to [bend left=25] ($(B)+(-\y,0)$);}
\draw[color=blue!60!black,rounded corners=5pt,line width=1pt] (0.5,0) -- (1.5,0.3)-- (3.6,0) -- (4.3,1.5)-- (4,3.2) -- (2.5,2.7) -- (0.2,3) -- (-0.2,2)-- cycle;
\end{tikzpicture}};
\node at (0,-2) {\begin{tikzpicture}[scale=0.5]
\draw [red] (1.1,0.75) -- (3.1,0.75);
\draw[color=blue!60!black,rounded corners=5pt,line width=1pt] (0.5,0) -- (1.5,0.15)-- (3.6,0) -- (4.3,0.75)-- (4,1.6) -- (2.5,1.35) -- (0.2,1.5) -- (-0.2,1)-- cycle;
\end{tikzpicture}};
\draw[->, color=blue!60!black] (0,-1) -- (0,-1.5);
\node at (0,-2.75) {\ding{193}};
\end{tikzpicture}
\end{array}
\end{array}
\]
where the contraction map $f$ is projective birational,  $\mathrm{R}f_*\scrO_{\scrX}=\scrR$, there is only one curve above the origin, $\scrX$ is smooth and the base $\scrR$ is complete local and Gorenstein.

\subsection{Noncommutative forms}\label{subsec:NCflopforms}
In classifying all $f$ with $\JRdim\Jac f\leq 1$, the main result \cite[A.13]{BW2} reduces us to one of the following three cases: either $\corank(f)\leq 1$, or $\corank(f)=2$ and either $\corank_{2}(f)=2,3$.  From this, a first coarse numerical subdivison is as follows:
\[
\begin{tabular}{ccccc}
\toprule
\textbf{Type} & $\corank(f)$ & $\corank_2(f)$  & $\corank_3(f)$ & $\corank_4(f)$ \\
\midrule
A & $\leq 1$ & $\leq 1$ & $\leq1$& $\leq1$\\
D & 2 & 2 & 2& $\leq2$ \\
$\mathrm{E}_6$ & 2 & 3 & 4& 4\\
\bottomrule
\end{tabular}
\]
To unpack this in one case, if $\corank(f)=2$ and $\corank_2(f)=3$, then up to isomorphism $f\in\mathbb{C}\llangle x,y\rrangle$ has zero linear and quadratic parts, and the derivatives of the cubic part are linearly dependent: that is, $f\cong x^3 +$ terms of order $\ge4$.

However, as stated the above is not quite the full picture. While these numerical conditions on a potential $f$ do happen to imply $\JRdim\Jac f\le1$ (but, only after classification), there are other families with $\JRdim\Jac f=0$ that are currently defined by more than just numerical coranks, using also analytic data about slices by central elements.
These other families are denoted $\mathrm{E}_{7}, \mathrm{E}_{8,5}$ and $\mathrm{E}_{8,6}$,
and \cite{BW3} explains both the motivation and the technical details.

Given this, we state the following slightly tersely, ignoring  the exceptional cases.
\begin{thm}[{\cite{BW2,BW3}}]\label{thm:mainNCformflop}
Situation \ding{192} (flopping).   Then $\mathrm{A}_{\mathrm{con}}\cong\Jac f$, where 
$f$ is either
\[
\begin{tabular}{p{1cm}p{3cm}p{4cm}}
$\mathrm{A}_n$& $x^{2}+y^{n}$ & $n\geq 2$\\[4pt]
$\mathrm{D}_{n,m}$ &$xy^{2}+x^{2m-1}+ x^{2n}$ & $n,m\geq 2$\\
$\mathrm{D}_{n,\infty}$ &$xy^{2}+x^{2n}$ & $n\geq2$\\[4pt]

$\mathrm{E}_{6,n}$&$x^{3}+xy^{3}+ p(y)$ & $p\in\mathbb{C}[t], \mathrm{ord}(p)=n\geq 4$ 
\end{tabular}
\]
or a normal form of Type $\mathrm{E}_{7}$, $\mathrm{E}_{8,5}$ or $\mathrm{E}_{8,6}$.
(These other Type~$\mathrm{E}$ forms also have moduli, but not the discrete parameter~$n$.)
\end{thm}

As stated, \ref{thm:mainNCformflop} is strictly speaking not a list of \emph{normal} forms, since different $p(y)$ in $\mathrm{E}_6$ may give isomorphic Jacobi algebras.  As in Arnold's original work \cite{Arnold}, the existence of moduli makes normal forms less pressing; the more important point is to cover all cases, which \ref{thm:mainNCformflop} does.  Moreover, there are indeed nontrivial moduli in $\mathrm{E}_6$, as we show in \cite{BW3} that two general members of the family $x^3 +xy^3 + y^4 + ay^5 + by^6$ with $a,b\in\mathbb{C}$ are not isomorphic.

In contrast to \ref{thm:mainNCformflop}, the following result does not need to be stated tersely: it turns out, contrary to our own initial expectations, that contraction algebras of crepant divisor-to-curve contractions are very rare, and indeed all can be defined numerically.  
\begin{thm}[{\cite{BW2,BW3}}]\label{thm:mainNCformdiv}
Situation \ding{193} (div$\to$curve).  Then $\mathrm{A}_{\mathrm{con}}\cong\Jac f$, where 
$f$ is either
\[
\begin{tabular}{p{1cm}p{3cm}p{4cm}}
$\mathrm{A}_\infty$& $x^{2}$ & \\[4pt]
$\mathrm{D}_{\infty,m}$ &$xy^{2}+ x^{2m-1}$ & $m\geq 2$\\
$\mathrm{D}_{\infty,\infty}$ &$xy^{2}$ &\\[4pt]

$\mathrm{E}_{6,\infty}$&$x^{3}+xy^{3}$ 
\end{tabular}
\]
\end{thm}
It is worth pausing to reflect that there are very few examples, namely just one countable family and three exceptional cases.   It seems surprising that $\JRdim\Jac f=1$ algebras are significantly more rare than those satisfying $\JRdim\Jac f=0$.  The second remarkable aspect of \ref{thm:mainNCformdiv} is that every case can be viewed as the $n\to\infty$ limit of a family in \ref{thm:mainNCformflop}.  This fact is a consequence of the classification, not an input, which we think of as trying to make sense in cold-blooded algebra of the vague intuition that every crepant divisor-to-curve contraction is an infinite limit of flops.

\begin{remark}\label{rem:Hao}
In addition to the results presented here, the upcoming paper of Zhang \cite{Hao} considers the noncommutative singularity theory on the following quiver, where the single loop at each vertex is optional. 
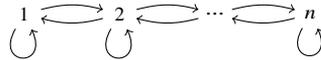
\[
\begin{tikzpicture}[scale=1.25,bend angle=15, looseness=1]
\node (a) at (-1,0) [pvertex] {$1$};
\node (b) at (0,0) [pvertex] {$2$};
\node at (0,0.6) [pvertex] {};
\node (c) at (1,0) [pvertex] {$\phantom{2}$};
\node at (1,0) {$\scriptstyle \hdots$};
\node (n) at (2,0) [pvertex] {$n$};
\draw[->,bend left] (b) to (a);
\draw[->,bend left] (a) to (b);
\draw[->,bend left] (c) to (b);
\draw[->,bend left] (b) to (c);
\draw[->,bend left] (n) to (c);
\draw[->,bend left] (c) to (n);
\draw[<-]  (a) edge [in=-120,out=-65,loop,looseness=7]  (a);
\draw[<-]  (b) edge [in=-120,out=-65,loop,looseness=7]  (b);
\draw[<-]  (n) edge [in=-120,out=-65,loop,looseness=7]  (n);
\end{tikzpicture}
\]
Zhang \cite{Hao} gives various numerical and structural characterisations of `Type $\mathrm{A}$' potentials on this quiver, and amongst other things proves the realisation conjecture (\ref{conj:georeal} below) for this intrinsic class of potentials.
\end{remark}

\subsection{Structural Geometric Consequences}\label{subsec:geocor}

Due to \ref{thm:DWconj}, the normal forms in the previous subsection have immediate consequences, once we are able to build the geometry which realises each such potential.  The overarching idea to do this elegantly is to use \cite{BW2b}, which was itself heavily inspired by the string theory literature.

\begin{thm}[\cite{Pagoda}, \cite{BW3} ongoing]
Smooth $3$-fold single-curve flops $\scrX\to\mathrm{Spec}\,\scrR$ are classified.  The classification data can be stated in various forms: by the equation defining~$\scrR$, by describing the glue defining $\scrX$, or by the slice of the universal flop.
\end{thm}
The precise statement of the theorem above is work in progress -- noting that Reid's Pagoda \cite{Pagoda} already handles the Type~A cases -- in part because some pieces of the data are much harder to describe than others.  The equation defining $\scrR$ together with its class group is perhaps the most compelling description, but that seems to be the hardest of all. All the same, from the perspective as users of the data, having multiple descriptions of each flop is the most beneficial outcome; for a commutative analogue, compare \cite{fifteen}, either as motivation or provocation.  

The various data that can be used to finally classify flops does serve to illustrate, above all, the elegance of the noncommutative forms in \S\ref{subsec:NCflopforms}. To give a feeling of this, from $\mathrm{E}_6$ onwards the equations defining the base of the flop become vast and unwieldy; as one example, the following is an example of the base of an $\mathrm{E}_8$ flop:
\tikzstyle{xxd} = [xshift=-1.7ex,yshift=0.8ex]
\tikzstyle{xx} = [xshift=-1.2ex,yshift=0.8ex]
\tikzstyle{xxs} = [xshift=-0.7ex,yshift=0.5ex]
\begin{align*}
\Spec\scrR = \left(\ 2041x^{12}\tikz[xxd] \coordinate (A1);\right. -107109x^{11}\tikz[xxd] \coordinate (A2); +{1390131}x^{10}\tikz[xxd] \coordinate (A3);  -{6130}x^9\tikz[xx] \coordinate (A4);y  -{8432}x^9\tikz[xx] \coordinate (A5); \qquad\quad \\ 
{}+{157127}x^8\tikz[xx] \coordinate (A6);y
- 17x^8\tikz[xx] \coordinate (A7);z  +{45}x^8\tikz[xx] \coordinate (A8);  
 -{654}x^7\tikz[xx] \coordinate (A9);y^2  -{340}x^7\tikz[xx] \coordinate (A10);y + {853}x^7\tikz[xx] \coordinate (A11);z \qquad\quad \\ 
{} -x^7\tikz[xx] \coordinate (A12);+  {21842}x^6\tikz[xx] \coordinate (A13);y^2 - {6}x^6\tikz[xx] \coordinate (A14);y + {5905}x^6\tikz[xx] \coordinate (A15);z + {7}{157}x^6\tikz[xx] \coordinate (A16);t -   {3}x^5\tikz[xx] \coordinate (A17);y^3 \qquad\quad \\ 
{} - {66}x^5\tikz[xx] \coordinate (A18);y^2 - {93}x^5\tikz[xx] \coordinate (A19);yz - {15}x^5\tikz[xx] \coordinate (A20);z - 
 {2358}x^5\tikz[xx] \coordinate (A21);t + {1090}x^4\tikz[xx] \coordinate (A22);y^3 + x^4\tikz[xx] \coordinate (A23);y^2 \qquad\quad \\ 
{} + {355}x^4\tikz[xx] \coordinate (A24);yz  - {9}x^4\tikz[xx] \coordinate (A25);z^2 + {7}x^4\tikz[xx] \coordinate (A26);t - x^3\tikz[xx] \coordinate (A27);y^3 - {2}x^3\tikz[xx] \coordinate (A28);y^2z -    {133}x^3\tikz[xx] \coordinate (A29);yt \qquad\quad \\
{}+ {15}x^3\tikz[xx] \coordinate (A30);z^2 + y^5 + {41}x^2\tikz[xx] \coordinate (A31);y^4 + {34}x^2\tikz[xx] \coordinate (A32);y^2z
 + {6}x^2\tikz[xx] \coordinate (A33);z^2 - x^2\tikz[xx] \coordinate (A34);zt \qquad\quad \\ 
 {} + x\tikz[xxs] \coordinate (A35);y^3z
  \left.- {15}x\tikz[xxs] \coordinate (A36);y^2t -  {5}x\tikz[xxs] \coordinate (A37);zt + z^3 + t^2=0 \ \right).
\end{align*}
The particular values of at least some of the coefficients in that mighty polynomial are important, but without significant work it is difficult to say which.
The polynomial would, all things considered, be hard to guess without the noncommutative normal form guiding its construction. It is at least reassuring to see that slicing the above by $x=0$ gives the equation $t^2 + z^3 + y^5=0$, which is the $E_8$ normal form in \ref{thm:ADE}.

\subsection{Numerical geometric consequences}\label{subsec:curveinv}
There are various numerical consequences of classification in the previous subsection. Fix a flopping contraction  $f\colon\scrX\to\mathrm{Spec}\,\scrR$ with exceptional curve $C$.  Then given any $\upbeta \in\mathbb{Z}[C]$, so just an integer, Katz \cite{KatzGV} defines an associated Gopakumar--Vafa (GV) invariant 
\[
n_{\upbeta}=n_{\upbeta}(\scrX) \in \mathbb{Z}_{\geq 0}
\]
by considering a one-parameter deformation
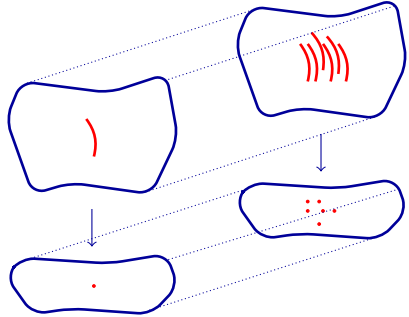
\[
\begin{array}{c}
\begin{tikzpicture}
\node at (0,0)
{$\begin{tikzpicture}
\node at (0,0) {\begin{tikzpicture}[scale=0.5]
\coordinate (T) at (1.9,2);
\coordinate (B) at (2.1,1);
\draw[red,line width=1pt] (T) to [bend left=25] (B);
\draw[color=blue!60!black,rounded corners=5pt,line width=1pt] (0.5,0) -- (1.5,0.3)-- (3.6,0) -- (4.3,1.5)-- (4,3.2) -- (2.5,2.7) -- (0.2,3) -- (-0.2,2)-- cycle;
\end{tikzpicture}};
\node at (0,-2) {\begin{tikzpicture}[scale=0.5]
\filldraw [red] (2.1,0.75) circle (1pt);
\draw[color=blue!60!black,rounded corners=5pt,line width=1pt] (0.5,0) -- (1.5,0.15)-- (3.6,0) -- (4.3,0.75)-- (4,1.6) -- (2.5,1.35) -- (0.2,1.5) -- (-0.2,1)-- cycle;
\end{tikzpicture}};
\draw[->, color=blue!60!black] (0,-1) -- (0,-1.5);
\end{tikzpicture}$};

\node at (3,1) {
$\begin{tikzpicture}
\node at (0,0) {\begin{tikzpicture}[scale=0.5]

\draw[red,line width=1pt] (1.5,2) to [bend left=25] (1.7,1);
\draw[red,line width=1pt] (1.7,2) to [bend left=25] (1.9,1);
\draw[red,line width=1pt] (1.8,2.3) to [bend left=25] (2.1,1.3);
\draw[red,line width=1pt] (2.1,2) to [bend left=25] (2.3,1);
\draw[red,line width=1pt] (2.2,2.2) to [bend left=25] (2.5,1.2);
\draw[red,line width=1pt] (2.5,2) to [bend left=25] (2.7,1);

\draw[color=blue!60!black,rounded corners=5pt,line width=1pt] (0.5,0) -- (1.5,0.3)-- (3.6,0) -- (4.3,1.5)-- (4,3.2) -- (2.5,2.7) -- (0.2,3) -- (-0.2,2)-- cycle;
\end{tikzpicture}};
\node at (0,-2) {\begin{tikzpicture}[scale=0.5]
\filldraw [red] (1.7,0.75) circle (1pt);
\filldraw [red] (2.1,0.75) circle (1pt);
\filldraw [red] (2.4,0.75) circle (1pt);
\filldraw [red] (1.7,1) circle (1pt);
\filldraw [red] (2,1) circle (1pt);
\filldraw [red] (2,0.4) circle (1pt);
\draw[color=blue!60!black,rounded corners=5pt,line width=1pt] (0.5,0) -- (1.5,0.15)-- (3.6,0) -- (4.3,0.75)-- (4,1.6) -- (2.5,1.35) -- (0.2,1.5) -- (-0.2,1)-- cycle;
\end{tikzpicture}};
\draw[->, color=blue!60!black] (0,-1) -- (0,-1.5);
\end{tikzpicture}$};
\draw[color=blue!60!black,densely dotted] (-0.85,-0.875) --(2.05,0.075);
\draw[color=blue!60!black,densely dotted] (1,-0.875) --(4.05,0.075);
\draw[color=blue!60!black,densely dotted] (0.75,-1.55) --(3.8,-0.55);

\draw[color=blue!60!black,densely dotted] (-0.85,1.475) --(2.1,2.45);
\draw[color=blue!60!black,densely dotted] (1,1.525) --(3.95,2.5);
\draw[color=blue!60!black,densely dotted] (0.7,0.025) --(3.7,1.025);

\end{tikzpicture}
\end{array}
\]
 such that
\begin{enumerate}
\item The central fibre $g_0\colon X_0\to Y_0$ is isomorphic to $f\colon\scrX\to\mathrm{Spec}\,\scrR$.
\item All other fibres $g_t\colon X_t\to Y_t$ are flopping contractions whose exceptional locus is a disjoint union of $(-1,-1)$-curves.
\end{enumerate}
The idea is then to obtain numerical information by counting the curves in a neighbouring fibre. 
Regarding the original flopping curve $C$ of $f$ as a curve in the fourfold $\mathcal{X}$ swept out by the one-parameter deformation, then the GV invariant $n_j$ is defined to be the number of $g_t$-exceptional $(-1,-1)$-curves $C'$ with curve class $j[C]$. In other words, $n_j$ is the number of $(-1,-1)$-curves $C'$ such that for every line bundle $\scrL$ on $\mathcal{X}$, 
\[
\deg (\scrL|_{C'})=j \deg(\scrL|_C).
\] 
Almost all the $n_j$ are zero. In fact, at most $(n_1,\hdots,n_6)$ are non-zero.

The following is a numerical consequence of the classification of flops via Jacobi algebras in the previous subsections.
\begin{thm}
The GV invariants are known for all smooth irreducible $3$-fold flops.  
\end{thm}
Full details will again appear in \cite{BW3}.  To outline  one particularly remarkable phenomenon, \emph{all} $\mathrm{E}_7$ flops (equivalently, the so-called length four flops) have the same GV invariants, namely $(6,5,2,1)$.  No other values are possible. Similarly, all the Gromov--Witten invariants and all the Donaldson--Thomas invariants are the same for any $\mathrm{E}_7$ flop. Thus from a numerical perspective, all $\mathrm{E}_7$ flops are identical: it is the extra information of the isomorphism class of the contraction algebra -- that is, the class of the noncommutative potential $f\in\mathbb{C}\llangle x, y\rrangle$ -- that distinguishes them.

\subsection{Conjectures}\label{subsec:conjectures}

We say that $f$ is \emph{geometric} if it arises from geometry, namely if $\Jac f$ is isomorphic to the contraction algebra of some  $\scrX\to\mathrm{Spec}\,\scrR$ in~\ref{deform curves} (either case \ding{192} or \ding{193}).

\begin{conj}[Geometric Realisation Conjecture]\label{conj:georeal}
Every $f$ whose Jacobi algebra satisfies $\dim\Jac f\leq 1$ is geometric.
\end{conj}

Tens of thousands of computer searches (with Magma \cite{magma} and Singular \cite{singular}) have failed to uncover counterexamples to \ref{conj:georeal}, but the methods of \cite{BW2, BW3} do not yet seem strong enough to prove it. At this stage, it seems that a new ingredient in noncommutative singularity theory will be needed, if indeed \ref{conj:georeal} is true. 

\medskip
Now let $C\cong\mathbb{P}^1\subseteq X$, where $X$ is a Calabi--Yau $3$-fold.  Write $\Lambda_{\mathrm{def}}$ for the noncommutative deformation algebra associated to $C$, obtained via noncommutative deformation theory. It is conjectured in \cite{BW2b} that in this 3-fold setting, numerical  properties of $\Lambda_{\mathrm{def}}$ are sufficient to detect whether or not $C$ contracts.

\begin{conj}[Artin Contractibility for $3$-folds]\label{conj:artin}
In a suitably local neighbourhood of $C$ in $X$, the following are equivalent.
\begin{enumerate}
\item $C$ contracts to a point and is an isomorphism elsewhere.
\item $\dim_\mathbb{C}\Lambda_{\mathrm{def}}<\infty$.
\end{enumerate}
\end{conj}
The statement of \ref{conj:artin} is for a single curve, but one of the advantages of the noncommutative deformation theory technology is that it finally permits clean statements for multi-curve configurations; for details see \cite{BW2b}.

It may be the case that the condition $\mathrm{H}^1(\hat{X},\scrO_{\hat{X}})=0$ needs to be added to the assumptions in the conjecture, since this would then more closely match the classical Artin statement for surfaces.  However, just as for \ref{conj:georeal} above, tens of thousands of computer searches suggest that the condition $\dim_\mathbb{C}\Lambda_{\mathrm{def}}<\infty$ is so restrictive that the rational neighbourhood assumption may in fact come for free.

\medskip
To end, it seems fitting to note the limitations of our view, or, perhaps better, the enormous potential in what has not been covered.
The discussion here has considered what we might term noncommutative function germs, that is, elements of $\mathbb{C}\llangle x_1,\hdots,x_d\rrangle$ or more generally elements of complete local path algebras on quivers. That is based largely on our admiraction and profitable exploitation of a few of the ideas of Arnold's seminal paper \cite{Arnold}. Conventional commutative singularity theory is of course a much larger subject, exemplified by the vast literature, such as \cite{ArnoldBook,mondbook,singhandbook} to cite a few thousand pages. We have not mentioned notions of stability, $k$\nobreakdash-determinacy, adjacency, versality, more general map germs, vector fields, or a host of other very natural notions. These may or may not have profitable extensions into the noncommutative realm. But for now we simply do not know, in part due to our technical limitations, but also since we do not yet know how we would be able to apply them.


\begin{ack}
Truth be told, it was during an otherwise uneventful late afternoon walk in September 2020 between the two authors and a young and crying baby, just past the junction of the footpath that links Palmerston Road to Beechwood Avenue in the suburbs of Coventry, where the idea of a noncommutative singularity theory was born.  Life is nothing, if not glamorous.  However, in that moment, absolutely everything we had been thinking about for almost a decade vastly simplified, and the road forward became crystal clear.

Then, but also much more broadly, maths is fundamentally a collaborative endeavour, and it would have been impossible to be delivering this article without the generous input and insights of many people over the years: from collaborators 
Jenny August,
Will Donovan,
Wahei Hara, 
Yuki Hirano, 
Osamu Iyama,
Navid Nabijou,  
Ivan Smith;
to PhD students 
Joe Karmazyn, Noah White, Kellan Steele, 
Okke van Garderen, Sarah Kelleher, 
Caroline Namanya, 
Brian Makonzi,  
Sam Lewis, 
Zhang Hao, 
Charlotte Bartram,
Marina Godinho, Parth Shimpi; 
and to postdocs 
Alice Rizzardo, 
Agnieszka Bodzenta, 
Theo Raedschelders, 
Diletta Martinelli, 
Ana Garcia Elsener, 
Franco Rota, 
Matt Pressland, 
Emma Lepri, and 
Timothy De Deyn.
\end{ack}

\begin{funding}
This work was supported by EPSRC grant EP/R034826/1 and by the ERC Consolidator Grant 101001227 (MMiMMa).
\end{funding}


\end{document}